\documentclass[12pt]{article}
\usepackage{latexsym, amsmath, amssymb, amsfonts, amscd}
\usepackage{color}
\usepackage{graphics}
\usepackage[dvips]{epsfig}

\usepackage[pagebackref]{hyperref}

\newtheorem{theorem}{Theorem}
\newtheorem{proposition}{Proposition}

\newtheorem{corollary}{Corollary}

\newtheorem{definition}{Definition}

\newcommand{\calo}{{\cal O}}

\newcommand{\adL}{\mbox{\rm ad}_{\Lambda}}

\newcommand{\ad}{\mbox{\rm ad}}

\def\lcf{\lbrack\! \lbrack}
\def\rcf{\rbrack\! \rbrack}

\def\dbar{\overline\partial}

\def\dbar{\overline\partial}

\def\oomega{\overline\omega}

\newcommand{\lieg}{\mathfrak{g}}
\newcommand{\liet}{\mathfrak{t}}
\newcommand{\liec}{\mathfrak{c}}

\newcommand{\lbra}[2]{\lcf #1, #2 \rcf}

\newcommand{\bproof}{\noindent{\it Proof: }}
\newcommand{\eproof}{\hfill \qed \vspace{0.2in}}
\def\qed{\rule{2.3mm}{2.3mm}}

\begin{document}
\title{\bf Algebraic Structure of  \\ Holomorphic Poisson Cohomology \\ on Nilmanifolds}
\author{
Yat Sun  Poon\thanks{ Address:
    Department of Mathematics, University of California,
    Riverside, CA 92521, U.S.A.. E-mail: ypoon@ucr.edu.}
    \
    and
John Simanyi\thanks{ Address:
    Department of Mathematics, University of California,
    Riverside, CA 92521, U.S.A.. E-mail: simanyi@math.ucr.edu.}
    }
%\date{September 3, 2018}
\maketitle
\begin{abstract} It is proved that on nilmanifolds with abelian complex structure, there exists
a canonically constructed non-trivial holomorphic Poisson structure. We identify the necessary and sufficient
condition for its associated cohomology to be isomorphic to the cohomology associated to trivial (zero) holomorphic Poisson structure.
We also identify a sufficient condition for this isomorphism to be at the level of Gerstenhaber algebras.
\end{abstract}

\section{Introduction}

The study of deformation of real Poisson bracket and its related cohomology theory was initiated by
Lichnerowicz et.~al.~long ago \cite{FLS} \cite{Lich}. Deformation of Poisson brackets in complex algebraic category emerged 20 years ago 
 \cite{Polish}. There has been associated deformation theory within the category of complex analytic geometry 
\cite{B-Marci} \cite{Hitchin-holomorphic Poisson} \cite{Ran}.

In the past fifteen years, holomorphic Poisson structure is also conceived as an object in the realm of generalized geometry \cite{Marco} 
\cite{Hitchin-Generalized CY} \cite{Hitchin-Instanton} .
In this vein, one studies  deformation within the framework of Lie bi-algebroids \cite{LWX}.
There has been work along this line by various authors \cite{Goto1} \cite{Goto2} \cite{GPR}
\cite{Hong} \cite{Hong-Xu} \cite{Xu}
\cite{PS}. In this note, the authors continue to treat holomorphic Poisson structures within the realm of 
generalized complex structures,  and extend the work in \cite{PS}.

A holomorphic bi-vector field $\Lambda$ is Poisson if its self bracket is equal to zero: $\lbra{\Lambda}{\Lambda}=0$. 
When $\Lambda=0$, it is a trivial case and we get nothing but a complex manifolds. As the condition $\lbra{\Lambda}{\Lambda}=0$ is homogeneous,
when $\Lambda\neq 0$,  
for any complex number $t$, $t\Lambda$ is a holomorphic Poisson structure. When $t\neq 0$, they are all equivalent and could 
be considered as a deformation of the trivial holomorphic Poisson structure given by $t=0$. Associated to each non-trivial holomorphic Poisson structure, one has 
a holomorphic version of Lichnerowicz's differentials $\dbar_{t\Lambda}$ 
\cite{CGP} \cite{Xu} \cite{Lich}. They act on the space of sections of exterior bundles of $(1,0)$-vectors and $(0,1)$-forms.
\[
\dbar_{t\Lambda}: \wedge^k(T^{1,0}M\otimes T^{*(0,1)}M) \to \wedge^{k+1}(T^{1,0}M\otimes T^{*(0,1)}M),
\]
for any non-negative integer $k$. They form an elliptic complex, and yield cohomology spaces. $\dbar_{t\Lambda}$ is a perturbation of the classical operator 
$\dbar$. When $t=0$, the resulting cohomology spaces are 
\[
H^k_{0}=\oplus_{p+q=k}H^q(M, \Theta^p)
\]
where $\Theta^p$ is the sheaf of germs of sections of the $p$-th exterior power of the holomorphic tangent bundle of the 
underlying complex manifold $M$. 
For each $t\neq 0$, we have cohomology spaces $H^k_{{t\Lambda}}$. They are all isomorphic to $H^k_{{\Lambda}}$.

Similar to the classical Lichnerowicz differential, $\dbar_{\Lambda}=\dbar+\adL$ where 
$\adL$ is the Schouten bracket of the Poisson structure $\Lambda$ with tensorial objects. 
In the holomorphic setting, the cohomology $H^k_{{\Lambda}}$ is the hyper-cohomology of a spectral sequence of a bi-complex. 
$\adL$ is also the operator on the first page of the spectral sequence: 
\[
\adL=d_1^{p,q}: H^q(M, \Theta^p)\to H^q(M, \Theta^{p+1}).
\]

The general goal of this piece of work is to understand the relation between  the algebraic structures on $H^k_{{\Lambda}}$ and those on $H^k_{0}$. 
 Our focus is on nilmanifolds with abelian complex structures. 
As found in \cite{CFG} \cite{CFGU} \cite{Salamon} and explained in Section \ref{ACS} in this paper,  
on such a manifold there exists a global frame of $(1,0)$-forms $\{\omega^1, \dots, \omega^n\}$ such that for each $k$,
$d\omega^k$ is type-$(1,1)$ and it is in the exterior algebra generated by $\{\omega^1, \dots, \omega^{k-1}, \oomega^1, \dots, \oomega^{k-1}\}$. We call 
such a basis an ascending basis. Let $\{V_1, \dots, V_n\}$ be the dual basis. Its complex linear span is denoted by $\lieg^{1,0}$.
As to be explained in  Definition \ref{existence} below, 
we find that $\Lambda=V_n\wedge V_{n-1}$ is a canonically defined holomorphic Poisson structure.

Our first main result is to identify a necessary and sufficient condition for the map $d_1^{p, q}$ of the spectral sequence associated to 
this particular holomorphic Poisson structure to be identically zero. 

\

\noindent\textbf{Theorem } {\it Let $M$ be a nilmanifold with abelian complex structure.
Let $\Lambda$ be the canonical holomorphic Poisson structure associated to an ascending basis of the complex structure.
The map $d_1^{p,q}$ is identically zero for all $p, q$ if and only if there exists  a vector $X$ in $\lieg^{1,0}$ such that
$\adL \oomega^n=\dbar X$. Whenever the vector $X$ exists, we call it a potential vector.}

\

The above theorem appears as Theorem \ref{degeneracy} in Section \ref{canonical}. Additional effort is needed to prove that 
the degeneracy of the spectral sequence actually lead to a vector space isomorphism
\begin{equation}\label{mapping}
\phi: 
\oplus_{p+q=k}H^q(M, \Theta^p) \to H^k_{\Lambda}
\end{equation}
as stated in Theorem \ref{surjective}. 
The proof of these two results occupies much of Section \ref{canonical}. They generalize a result in \cite{PS}.

As explained in Proposition \ref{characterization}, 
given $\Lambda=V_n\wedge V_{n-1}$, a characterization of the potential vector $X$
is  the set of equations below.
\[
\iota_Xd\omega^k=0, \quad \mbox{ for all } \quad 1\leq k\leq n-1, \quad \mbox{ and } \quad 
\iota_Xd\omega^n=-\iota_{V_{n-1}}d\oomega^n. 
\]
This characterization also corrects an error in Proposition 9 of \cite{PS}. 

After the conditions on a vector space isomorphism is clarified, we complete an analysis in Theorem \ref{algebra iso}
 on whether the isomorphism $\phi$ in (\ref{mapping}) is an isormophism of algebras. We state the result below. 

\

\noindent\textbf{Theorem } {\it  Suppose that $\Lambda$ is the canonical holomorphic Poisson
structure associated to an ascending basis of an abelian
complex structure on a nilmanifold. If  $X$ is a potential vector, 
then the map $\phi$ is an isomorphism of Gerstenhaber algebras when $\ad_X\oomega^k=\lbra{X}{\oomega^k}=0$ for  all
$1\leq k\leq n-1$.}

\

The above result upgrades our past observation regarding isomorphism on the level of vector spaces. 
As an application of the above results, we prove Theorem \ref{key} in Section \ref{two-steps case} as stated below.

\

\noindent\textbf{Theorem } {\it 
On a two-step nilmanifold with abelian complex structure, if $d\omega^n$ is  a non-degenerate map then 
the canonical Poisson cohomology with respect to an
ascending basis is isomorphic to the one associated to the trivial (zero) Poisson structure as Gerstenhaber algebras. }

\

The above result recovers an ad hoc computation by the first author in his paper \cite{Poon} on Kodaira surfaces, 
and give the computation a theoretical understanding. 
It also extends and explains  various examples of 
two-steps nilmanifolds found in \cite{PS}. At the end of this article we provide abundant examples beyond two-step nilmanifolds
on which the conditions for Theorem \ref{algebra iso} are satisfied. We also found
examples for which the conclusion for Theorem \ref{degeneracy} does not hold. Together, they present a rich subject
for further investigation.

\section{Algebraic Structures on Poisson Cohomology}\label{ASPC}
Let $M$ be a manifold with an integrable complex structure $J$.
The complexified tangent bundle of $M$ splits into direct sum of the bundle of $(1,0)$-vectors $TM^{1,0}$ and the bundle of
$(0,1)$-vectors $TM^{0,1}$. We denote their dual bundles by $TM^{*(1,0)}$ and $TM^{*(0,1)}$ respectively.
Their $p$-th exterior bundles are denoted by
$TM^{p,0}$, $TM^{0,p}$, $TM^{*(p,0)}$ and $TM^{*(0,p)}$ respectively.
Define $L=TM^{1,0}\oplus TM^{*(0,1)}$. Via a natural pairing between forms and vectors, its complex conjugation $\overline{L}$ is complex linearly
isomorphic to its dual $L^*$. The pair $L$ and $\overline{L}$ forms a complex Lie bi-algebroid with respect to the Courant bracket \cite{Gau} \cite{LWX}.
It follows the theory of Lie bi-algebroid \cite{LWX} \cite{Mac} and there exists a natural differential
\begin{equation}
\dbar: C^\infty(M, L ) \to
C^\infty(M, \wedge^2L).
\end{equation}
When the Courant bracket is restricted to $L$, it is the Schouten bracket \cite{Mac},
also known as Fr\"olicher-Nijenhuis bracket in classical complex deformation theory.
It is extended to a differential of exterior algebras:
\begin{equation}
\dbar: C^\infty(M, \wedge^kL) \to
C^\infty(M, \wedge^{k+1}L).
\end{equation}
The differential $\dbar$ has the properties such that  when it is restricted to
$(0,1)$-forms, it is the classical $\dbar$-operator in complex manifold theory; meaning that it
is the $(0,2)$-component of the exterior differential \cite{GPR}.
Similarly, when the Lie algebroid differential is restricted to
$(1,0)$-vector fields, it is
the  Cauchy-Riemann operator as seen in \cite{Gau}. The differential $\dbar$, the bracket $\lbra{-}{-}$ and the exterior product $\wedge $ together equip the
space of sections of the exterior bundle $C^\infty(M, \wedge^\bullet L)$ with the structure of a differential Gerstenhaber algebra. In particular,
 if $a$ is a smooth section
of  $\wedge^{|a|}L$, where $|a|$ denotes the degree of $a$ in the exterior algebra,
and $b$ is a smooth section of $\wedge^{|b|}L$, then
\begin{eqnarray}
\dbar \lbra{a}{b} &=& \lbra{\dbar a}{ b}+(-1)^{|a|+1}\lbra{a}{\dbar b};\\
\dbar (a\wedge b) &=& (\dbar a)\wedge b+(-1)^{|a|}a\wedge (\dbar b),
\end{eqnarray}

Since $\dbar\circ\dbar=0$, one obtains the Dolbeault cohomology with coefficients in holomorphic
polyvector fields. Denoting the sheaf of germs of sections of
the $p^{\mathrm{th}}$ exterior power of the holomorphic tangent bundle by $\Theta^p$, we have
\[
H^\bullet_0=H^\bullet_{\dbar}(M, \wedge^\bullet L)\cong
\bigoplus_{p,q\geq 0}H^q(M, \Theta^p).
\]
In subsequent computations, when $p=0$, $\Theta^p$ represents the structure
sheaf $\mathcal{O}$ of the complex manifold $M$.

Due to the compatibility between $\dbar$ and
 the Schouten bracket $\lbra{-}{-}$, and the compatibility between $\dbar$
 and the exterior product $\wedge$
 as noted above,  the Schouten bracket and exterior product descend to
the cohomology space $H^\bullet(M, \Theta^\bullet)$. In other words,
the triple
\begin{equation}\label{model}
(\oplus_{p,q}H^q(M, \Theta^{p}), \lbra{-}{-}, \wedge)
\end{equation}
 forms a Gerstenhaber algebra.

 Suppose that $\Lambda$ is a holomorphic Poisson bivector field, then $\lbra{\Lambda}{\Lambda}=0$ and $\dbar\Lambda=0$.
 Consider the map
 \begin{equation}
 \adL:  C^\infty(M, \wedge^kL)\to C^\infty(M, \wedge^{k+1}L)
 \end{equation}
 defined by $\adL(\Phi)=\lbra{\Lambda}{\Phi}$ where $\Phi\in C^\infty(M, \wedge^kL)$.
 The identity $\lbra{\Lambda}{\Lambda}=0$ is translated into $\adL\circ\adL=0$. Since the complex structure
 is integrable, we have $\dbar\circ \dbar=0$. Therefore, the map
 \begin{equation}
 \dbar_{\Lambda}=\adL+\dbar: C^\infty(M, \wedge^kL)\to C^\infty(M, \wedge^{k+1}L)
 \end{equation}
 satisfies the condition $\dbar_{\Lambda}\circ \dbar_{\Lambda}=0$.
 It generates a complex, and hence a cohomology space. We denote it by $H^\bullet_{\Lambda}(M)$, and
 call it the (associated) Poisson cohomology.
 Furthermore, as $\Lambda$ is a degree-2 element in the differential Gerstenhaber algebra
 $(C^\infty(M, \wedge^\bullet L), \dbar, \lbra{-}{-}, \wedge )$, the operator
 $\adL$ satisfies the identities:
 \begin{eqnarray*}
\adL \lbra{a}{b} &=& \lbra{\adL a}{ b}+(-1)^{|a|+1}\lbra{a}{\adL b};\\
\adL (a\wedge b) &=& (\adL a)\wedge b+(-1)^{|a|}a\wedge (\adL b).
\end{eqnarray*}
So does $\dbar_{\Lambda}$.
\begin{eqnarray}
\dbar_{\Lambda} \lbra{a}{b} &=& \lbra{\dbar_{\Lambda} a}{ b}+(-1)^{|a|+1}\lbra{a}{\dbar_{\Lambda} b};\\
\dbar_{\Lambda} (a\wedge b) &=& (\dbar_{\Lambda} a)\wedge b+(-1)^{|a|}a\wedge (\dbar_{\Lambda} b).
\end{eqnarray}
It follows that $(C^\infty(M, \wedge^\bullet L), \dbar_{\Lambda}, \lbra{-}{-}, \wedge )$ is a differential Gerstenhaber algebra, and
$(H^\bullet_\Lambda(M),  \lbra{-}{-}, \wedge )$ inherits the structure of a Gerstenhaber algebra.
Observe that
\[
\wedge^kL=\oplus_{p+q=k}TM^{p,0}\otimes TM^{*(0,q)}
\]
and
\begin{eqnarray}
&\adL: C^\infty(M, TM^{p,0}\otimes TM^{*(0,q)})\to C^\infty(M, TM^{p+1,0}\otimes TM^{*(0,q)}), \label{adL on pq} & \\
&\dbar: C^\infty(M, TM^{p,0}\otimes TM^{*(0,q)}) \to C^\infty(M, TM^{p,0}\otimes TM^{*(0,q+1)}).&
\end{eqnarray}
Since
$\dbar\circ \dbar=0,$  $ \adL\circ \dbar+\dbar\circ\adL=0,$ and $\adL\circ\adL=0,$
the cohomology $H^\bullet_\Lambda(M)$ is given by the hypercohomology of the bi-complex $ C^\infty(M, TM^{\bullet,0}\otimes TM^{*(0,\bullet)})$.
The first page of related spectral sequence is given by the classical Dolbeault cohomology of $(0,q)$-forms with coefficients in the sheaf of germs
of holomorphic $p$-vector fields \cite{Voisin}.
\begin{equation}\label{E1}
E_1^{p,q}=H^q(M, \Theta^p).
\end{equation}
For $p=0$, $E_1^{0,q}=H^q(M, \calo)$, the space of holomorphic $(0,q)$-forms.
The map $d_1^{p,q}$ in the first page of the spectral sequence is the map $\adL$, which descends from 
(\ref{adL on pq}) to cohomology level. 
\begin{equation}
\adL: H^q(M, \Theta^p) \to H^q(M, \Theta^{p+1}).
\end{equation}

\section{Abelian Complex Structure on Nilmanifolds}\label{ACS}
After reviewing some fundamental observations on complex structures on nilmanifolds as
 given in \cite{Salamon}, we set up several technical observations
to facilitate the computation in subsequent sections.

A nilmanifold $M$ is a co-compact quotient of a simply connected nilpotent Lie group $G$.
Denote the Lie algebra of the group $G$ by $\lieg$. We assume that $\dim_R\lieg=2n$. Denote the dual space by
$\lieg^*$. One has the Chevalley-Eilenberg  differential $d$. If $\alpha$ is in $\lieg^*$, $X,Y$ are
in $\lieg$, then by definition
\begin{equation}\label{EC differential}
d\alpha (X,Y)=-\alpha (\lbra{X}{Y}).
\end{equation}
It is extended to the Chevelley-Eilenberg complex,
\[
0\to \lieg^* {\stackrel{d}\longrightarrow} \wedge^2\lieg^* {\stackrel{d}\longrightarrow}\cdots
{\stackrel{d}\longrightarrow}\wedge^{2n}\lieg^*\to 0.
\]
The descending central series of $\lieg$ is the chain of ideals defined inductively by $\lieg^0=\lieg$ and
$\lieg^{j}=[\lieg^{j-1}, \lieg]$ for all $j\geq 1$. By definition, the Lie algebra $\lieg$ is $s$-step
nilpotent if $\lieg^s=\{0\}$ and $\lieg^{s-1}\neq \{ 0\}$. In particular, the chain of ideals terminates,
\[
\lieg=\lieg^0\supseteq \lieg^1\supseteq\cdots \supseteq\lieg^{s-1}\supseteq\lieg^s=\{ 0\}.
\]
A real linear map $J:\lieg \to \lieg$ with the properties that $J\circ J=-$identity and for all
$X,Y\in \lieg$,
\begin{equation}\label{abelian cx}
\lbra{JX}{JY}=\lbra{X}{Y}
\end{equation}
is an abelian complex structure. Denote the spaces of $(1,0)$-vectors and $(0,1)$-vectors
respectively by $\lieg^{1,0}$ and $\lieg^{0,1}$.
The identity (\ref{abelian cx}) above implies that
$\lieg^{1,0}$ is an abelian complex algebra. Therefore, the only Lie bracket among elements of
$\lieg_C$ is of the form $\lbra{X}{\overline Y}$ and its complex conjugate, where $X\in \lieg^{1,0}$ and ${\overline Y}\in \lieg^{0,1}$.

Taking the adjoint of $J$, one obtains a map, also denoted by $J$, from
$\lieg^*$ to itself. Denote the space of  $(1,0)$-forms by
\[
\lieg^{*(1,0)}=\{\alpha-iJ\alpha: \alpha\in \lieg^*\}.
\]
Its complex conjugation is the space of  $(0,1)$-forms. It is
denoted by $\lieg^{*(0,1)}$.

Their $p$-th exterior products are respectively denoted by $\lieg^{p,0}$,
$\lieg^{0,p}$, $\lieg^{*(p,0)}$ and $\lieg^{*(0,p)}$. In particular, by taking invariant polyvector fields and invariant forms,
we have an inclusion map
\begin{equation}
\iota: \lieg^{p,0}\otimes \lieg^{*(0,q)}\to C^\infty(M, TM^{p,0}\otimes TM^{*(0,q)}).
\end{equation}

\begin{theorem}\label{algebra}
 {\rm \cite{CFP} } The inclusion map $\iota$ is a quasi-isomorphism. In other words, it induces an isomorphism of cohomology
\[
\iota: H^q(\lieg^{p,0})\cong H^q(M, \Theta^p).
\]
\end{theorem}

The Chevelley-Eilenberg complex is complex linearly extended.
When the complex structure is abelian as defined in (\ref{abelian cx}), then $d\alpha$ is a type-(1,1) form
when $\alpha$ is a type-(1,0) form. In the rest of this article, the following observation will play a
crucial role.

\begin{theorem}\label{Salamon basis} {\rm \cite{CFG} \cite{CFGU} \cite{Salamon}\ } 
If $\lieg$ is a nilpotent Lie algebra with an abelian complex structure,
then there exists a basis $\{\omega^1, \dots, \omega^n\}$ of $\lieg^{*(1,0)}$ such that
\[
d\omega^{j+1}\in I(\omega^1, \dots, \omega^j)\wedge I(\oomega^1, \dots, \oomega^j),
\]
where $I(\omega^1, \dots, \omega^j)$ denotes the ideal generated by $\{\omega^1, \dots, \omega^j\}$.
\end{theorem}
In other words, the structure constants $A_{k\ell}^j$ are such that for each $j$,
\begin{equation}\label{d omega j}
d\omega^j=\sum_{1\leq k,\ell\leq j-1}A^{j}_{k\ell}\omega^k\wedge\oomega^\ell.
\end{equation}
Taking complex conjugation and rearranging indices, we get
\begin{equation}\label{d oomega j}
d\oomega^j=-\sum_{1\leq k,\ell\leq j-1}{\overline A}^{j}_{\ell k}\omega^k\wedge\oomega^\ell.
\end{equation}
In particular,
\begin{equation}\label{d omega n}
d\omega^n=\sum_{1\leq k,\ell\leq n-1}A^{n}_{k\ell}\omega^k\wedge\oomega^\ell,
\quad
\mbox{ and } \quad
d\oomega^n=-\sum_{1\leq k,\ell\leq n-1}{\overline A}^{n}_{\ell k}\omega^k\wedge\oomega^\ell.
\end{equation}

As noted in \cite{Salamon}, one obtains the basis $\{\omega^1, \dots, \omega^n\}$ inductively starting
 with $d\omega^1=0$. Therefore,
we address this basis for $\lieg^{1,0}$ as an \it ascending basis. \rm Apparently, an ascending basis is not necessarily unique.
Let $\{V_1, \dots, V_n\}$ be the dual basis for $\lieg^{1,0}$. We will refer to it as the (dual) ascending basis
for $\lieg^{1,0}$.

We will also use the following notations,
\begin{eqnarray}
\liet^{1,0}=\mbox{span of }\{V_1, \dots, V_{n-1}\},
&
\liec^{1,0}=\mbox{span of }\{V_n \},\\
\liet^{*(1,0)}=\mbox{span of }\{\omega^1, \dots, \omega^{n-1}\},
&
\liec^{*(1,0)}=\mbox{span of }\{\omega^n \}.
\end{eqnarray}
Their conjugated counterparts $\liet^{0,1}, \liec^{0,1}, \liet^{*(0,1)}$ and
$\liec^{*(0,1)}$ are naturally defined. With these notations, we consider the contraction
of $d\omega^n$ and $d\oomega^n$ as linear maps.
\begin{equation}\label{linear map}
d\omega^n: \liet^{1,0}\to \liet^{*(0,1)},
\quad
d\oomega^n: \liet^{1,0}\to \liet^{*(0,1)}. 
\end{equation}
Moreover, with the given ordered bases, the matrix representation $A^n$
of $d\omega^n$ is $A^n_{k\ell}$ acting on row
vectors from the right. i.e. if $V=c^1V_1+\cdots+c^{n-1}V_{n-1}$, then
\[
\iota_Vd\omega^n=\sum_{1\leq k,\ell\leq n-1}c^kA^{n}_{k\ell}\oomega^\ell.
\]
For instance,
\begin{equation}\label{contract Vn-1}
\iota_{V_{n-1}}d\omega^n=\sum_{1\leq \ell\leq n-1}A^{n}_{n-1,\ell}\oomega^\ell.
\end{equation}
The matrix representation of $d\oomega^n$ is $-{\overline A}^T$.

\

From the structure equations (\ref{d omega j}) and (\ref{d oomega j}), one finds the expression of
structure equations in the basis $\{V_1, \dots, V_n\}\cup \{{\overline{V}}_1, \dots, {\overline{V}}_n\}$ for the complexified
Lie algebra $\lieg_{C}$. Namely,
\begin{equation}
\lbra{V_k}{{\overline V}_\ell}=-\sum_{k,\ell\leq j-1}A^j_{k\ell}V_j
+\sum_{k,\ell\leq j-1}{\overline A}^j_{\ell k}{\overline V}_j.
\end{equation}

It is now apparent that $V_n$ is in the center of $\lieg_C$.
Moreover, for each $\omega^k$ in the ascending basis and any $\overline {Y}$ in $\lieg^{0,1}$,
\[
(\dbar X)(\omega^k, \overline{Y})=X(\lbra{\overline{Y}}{\omega^k})=X(\iota_{\overline{Y}}d\omega^k)=-d\omega^k(X, \overline{Y})=-(\iota_Xd\omega^k)(\overline{Y}).
\]
It follows that
\begin{equation}\label{dbar X}
\dbar X=-\sum_{k=1}^n V_k\wedge (\iota_Xd\omega^k).
\end{equation}
As a result of (\ref{contract Vn-1}),
\begin{equation}
\dbar V_n=0, \quad \mbox{ and } \quad \dbar V_{n-1}=-V_n\wedge \iota_{V_{n-1}}d\omega^n.
\end{equation}
\begin{corollary} $\Lambda=V_n\wedge V_{n-1}$ is a holomorphic Poisson bivector.
\end{corollary}
\bproof $\Lambda$ is Poisson because the complex algebra $\lieg^{1,0}$ is abelian.
It is holomorphic because
\[
\dbar(V_n\wedge V_{n-1})=(\dbar V_n)\wedge V_{n-1}-V_n\wedge (\dbar V_{n-1})
=V_n\wedge V_n\wedge\left( \sum_{k=1}^{n-1}A^n_{n-1,k} \oomega^k\right)=0.
\]
\eproof

\begin{definition}\label{existence} Given an ascending basis $\{V_1, \dots, V_{n-1}, V_n\}$ for $\lieg^{1,0}$ 
on a nilmanifold with abelian complex
structure, we call $\Lambda=V_n\wedge V_{n-1}$ its canonical holomorphic Poisson structure, or simply a canonical Poisson structure of the
abelian complex structure.
\end{definition}

\section{Cohomology of Canonical Poisson Structures}\label{canonical}

Since the canonical Poisson structure is invariant, given the quasi-isomorphism in Theorem \ref{algebra},
one could reduce the spectral sequence computation with
$E_1$-terms given in (\ref{E1}) to a computation of invariant objects.
Namely, the $d_1$ map is reduced to
\begin{equation}\label{d1}
d_1^{p,q}=\adL: H^q(\lieg^{p,0})\to H^q(\lieg^{p+1,0}).
\end{equation}

The goal of this section is to prove the following theorem, which generalizes a key result in \cite{PS}.

\begin{theorem}\label{degeneracy} Let $M$ be a nilmanifold with abelian complex structure.
Let $\Lambda$ be the canonical holomorphic Poisson structure associated to an ascending basis of the complex structure.
The map $d_1^{p,q}$ is identically zero for all $p, q$ if and only if there exists  a vector $X$ in $\lieg^{1,0}$ such that
$\adL \oomega^n=\dbar X$.
\end{theorem}

We call such $X$, if exists, a potential vector for the canonical holomorphic Poisson structure. The proof below is similar to,
but simplifies an argument in \cite{PS} due to the adoption of an ascending basis.

\

\bproof
One direction of the computation is straightforward. Assume that
the map $d_1^{p,q}$ is identically zero for all $p, q$. When $(p, q)=(0, 1)$, the space $H^1(\lieg^{0,0})$ is precisely the space of invariant
closed $(0,1)$-forms. Given that the complex structure is abelian, it is apparent that every element in an
ascending basis
$\{\oomega^1, \dots, \oomega^n\}$ is $\dbar$-closed. Furthermore, given that it is an ascending basis,
\begin{equation}
\iota_{V_n}d\oomega^j=0 {\mbox{ for all }}  1\leq j\leq n, \quad \mbox{ and } \quad \iota_{V_{n-1}}d\oomega^k=0 {\mbox{ for all }} 1\leq k \leq n-1.
\end{equation}
Therefore, the only possible non-trivial term in the computation of $d_1^{0,1}$ is contributed by
\[
\adL\oomega^n=\lbra{V_n\wedge V_{n-1} }{\oomega^n}=V_n\wedge \iota_{V_{n-1}} d\oomega^n.
\]
It is equal to zero in $H^1(\lieg^{1,0})$ as a cohomology class if and only if there exists an element $X$ in $\lieg^{1,0}$ such that
\begin{equation}\label{dbar potential}
\adL(\oomega^n)=\dbar X.
\end{equation}

Next, we prove that the existence of a potential vector is a sufficient condition for the spectral sequence to degenerate.
Suppose $\Gamma^{p,q}$ is an element
in $\lieg^{p,0}\otimes\lieg^{*(0,q)}$ representing a
class in $H^q(\lieg^{p,0})$. Due to the direct sum decomposition
$\lieg^{*(0,1)}=\liec^{*(0,1)}\oplus \liet^{*(0,1)}$, we have
\begin{equation}\label{ct decomposition}
\lieg^{p,0}\otimes\lieg^{*(0,q)} = \lieg^{p,0}\otimes \wedge^q(\liec^{*(0,1)}\oplus \liet^{*(0,1)})
= \lieg^{p,0}\otimes (\liec^{*(0,1)}\otimes \liet^{*(0,q-1)} \oplus \liet^{*(0,q)}).
\end{equation}
Therefore,
\begin{equation}\label{pq decomposition}
\Gamma^{p,q}=\oomega^n\wedge \sum_{a}({\Pi}_a\wedge{\overline\Phi}_a)+\sum_b{\Upsilon}_b\wedge {\overline\Psi}_b,
\end{equation}
where $\{{\overline\Phi}_a, a=1, 2, \dots \}$ and $\{{\overline\Psi}_b, b=1, 2, \dots \}$ are generated by the exterior products of
$\{\oomega^1, \dots \oomega^n\}$, and they form bases for  $\liet^{*(0,q-1)}$ and  $\liet^{*(0,q)}$ respectively.
The terms $\Pi_a$ and $\Upsilon_b$ are elements in $\lieg^{p,0}$.

Since $\lbra{V_n}{\oomega^k}=\iota_{V_n}d\oomega^k=0$ for all $k$
and $\lbra{V_{n-1}}{\oomega^k}=\iota_{V_{n-1}}d\oomega^k=0$ for all $k\leq n-1$,
\begin{eqnarray*}
\adL(\Gamma^{p,q}) &=&\adL( \oomega^n\wedge \sum_{a}({\Pi}_\alpha\wedge{\overline\Phi}_a))
=\lbra{V_n\wedge V_{n-1}}{\oomega^n\wedge \sum_{a}({\Pi}_a\wedge{\overline\Phi}_a)}\\
&=&V_n\wedge\lbra { V_{n-1}}{\oomega^n\wedge \sum_{a}({\Pi}_\alpha\wedge{\overline\Phi}_a)}
=V_n\wedge\lbra { V_{n-1}}{\oomega^n}\wedge \sum_{a}({\Pi}_a\wedge{\overline\Phi}_a)\\
&=&\lbra{\Lambda}{\oomega^n}\wedge\sum_{a}({\Pi}_a\wedge{\overline\Phi}_a)
=(\dbar X)\wedge\sum_{a}({\Pi}_a\wedge{\overline\Phi}_a).
\end{eqnarray*}
On the other hand, since $\Gamma^{p,q}$ is $\dbar$-closed and all $(0, 1)$-forms of an ascending basis are also $\dbar$-closed,
\begin{eqnarray*}
\dbar \Gamma^{p,q} &=& \dbar( \oomega^n\wedge \sum_{a}({\Pi}_a\wedge{\overline\Phi}_a))
+\sum_b \dbar({\Upsilon}_b\wedge {\overline\Psi}_b)\\
&=&-\oomega^n\wedge \sum_{a}((\dbar{\Pi}_a)\wedge{\overline\Phi}_a)
+\sum_b (\dbar{\Upsilon}_b)\wedge {\overline\Psi}_b.
\end{eqnarray*}
In view of the nature of an ascending basis, every $(\dbar{\Upsilon}_b)\wedge {\overline\Psi}_b$ is an element in
$ \lieg^{p,0}\otimes \liet^{*(0,q+1)}$ while $\oomega^n\wedge \sum_{a}((\dbar{\Pi}_a)\wedge{\overline\Phi}_a)$  is in
$ \lieg^{p,0}\otimes \liec^{*(0,1)}\otimes \liet^{*(0,q)}.$ Therefore, $\dbar\Gamma^{p,q}=0$ if and only if
\[
\oomega^n\wedge \sum_{a}((\dbar{\Pi}_a)\wedge{\overline\Phi}_a)=0,
\quad {\mbox{ and  } } \quad
\sum_b (\dbar{\Upsilon}_b)\wedge {\overline\Psi}_b=0.
\]
In particular,
\[
\dbar\sum_{a}({\Pi}_a)\wedge{\overline\Phi}_a)=\sum_{a}((\dbar{\Pi}_a)\wedge{\overline\Phi}_a)=0.
\]
 It follows from the previous paragraph that
\begin{equation}
\adL(\Gamma^{p,q})=(\dbar X)\wedge\sum_{a}({\Pi}_a\wedge{\overline\Phi}_a)=
\dbar \left(X\wedge\sum_{a}({\Pi}_a\wedge{\overline\Phi}_a)\right).
\end{equation}
In terms of spectral sequence, $d_1^{p,q}(\Gamma^{p,q})=\adL(\Gamma^{p,q})\equiv 0$ as a cohomology class in $H^{q+1}(\lieg^{p,0})$,
concluding the proof of Theorem \ref{degeneracy}.
\eproof

In subsequent computation, given any element $\Gamma^{p,q}$ in $\lieg^{p,0}\otimes\lieg^{*(0,q)}$, we represent its decomposition
in (\ref{ct decomposition}) and (\ref{pq decomposition}) by
\begin{equation}\label{a n b}
\Gamma^{p,q}=\oomega^n\wedge \alpha^{p, q-1} +\beta^{p,q}.
\end{equation}
The proof above means that whenever $\dbar\Gamma^{p,q}=0$, then $\dbar\alpha^{p, q-1}=0$ and $\dbar\beta^{p,q}=0$, and
\begin{equation}
\adL(\Gamma^{p,q})=\dbar(X\wedge \alpha^{p, q-1})
\end{equation}
when $\adL\oomega^n=\dbar X$.

\begin{theorem}\label{surjective} When $\adL(\oomega^n)$ is $\dbar$-exact, there is a natural isomorphism
\[
\phi: \oplus_{p+q=k}H^q(\lieg^{p,0}) \to H^k_{\Lambda}.
\]
\end{theorem}
\bproof
Let
\begin{equation}
B^{p,q}=\lieg^{p,0}\otimes \lieg^{*(0,q)}.
\end{equation}
As noted in \cite{PS}, the bi-grading yields a natural filtration of the cohomology $H^n_\Lambda$:
\begin{equation}
H^k_\Lambda=F^0H^k\supset F^1H^k \supset \cdots \supset F^\ell H^k \supset  \cdots \supset F^kH^k \supset F^{k+1}H^k=\{0\}
\end{equation}
such that
\[
H^k_\Lambda\cong \oplus_{m=0}^k\frac{F^mH^k}{F^{m+1}H^k}.
\]
There is a natural injective map from
\[
\psi: \frac{F^mH^k}{F^{m+1}H^k}\to H^{k-m}(\lieg^{m,0})
\]
defined as below.
A non-trivial class in $\frac{F^mH^k}{F^{m+1}H^k}$ is presented by
\[
\Gamma=\sum_{m\leq p\leq k,  p+q=k}\Gamma^{p, q}=\Gamma^{k,0}+\Gamma^{k-1,1}+\cdots +\Gamma^{m+1, k-m-1}+\Gamma^{m, k-m}
\]
with $\Gamma^{m, k-m}\neq 0$ such that  $\dbar_\Lambda\Gamma=0$.
Recall that $\dbar_\Lambda=\adL+\dbar$, when we group $\dbar_\Lambda\Gamma$ by bi-degrees, we find
\begin{eqnarray*}
&& \dbar_\Lambda\Gamma\\
&=&\adL\Gamma^{k,0}+(\dbar\Gamma^{k,0}+\adL\Gamma^{k-1, 1})+
\cdots+(\dbar\Gamma^{m+1, k-m+1}+\adL\Gamma^{m, k-m})+\dbar\Gamma^{m,k-m}.
\end{eqnarray*}
Therefore, $\dbar\Gamma^{m, k-m}=0$, and hence 
$\Gamma^{m, k-m}$ represents a class in $H^{k-m}(\lieg^{m,0})$.
One could verify that it represents a well defined image of $\psi(\Gamma)$.

To prove that the map $\psi$ is surjective, recall the decomposition of $\Gamma^{p,q}$ as 
in (\ref{pq decomposition}), or equivalently (\ref{a n b}) when $q\geq 1$.
 \begin{eqnarray*}
&& \dbar_{\Lambda}(\Gamma^{p,q}- X \wedge \alpha^{p, q-1})\\
&=&\adL (\Gamma^{p,q}-X\wedge\alpha^{p, q-1})+\dbar (\Gamma^{p,q}-X\wedge\alpha^{p, q-1})\\
&=&\adL \Gamma^{p,q}-\dbar(X\wedge\alpha^{p, q-1})-\adL(X\wedge\alpha^{p, q-1})+\dbar \Gamma^{p,q} \\
&=& -\adL(X\wedge\alpha^{p, q-1}).
 \end{eqnarray*}
 Note that $\alpha^{p,q-1}$ is in $\lieg^{p,0}\otimes \liet^{*(0,q-1)}$ and the complex structure is abelian. The only non-trivial outcome in
 \[
\adL(X\wedge \alpha^{p, q-1})=\lbra{V_n\wedge V_{n-1} }{X\wedge \alpha^{p, q-1}}
\]
is contributed by $\lbra{V_n}{\oomega^k}$ and $\lbra{V_{n-1}}{\oomega^k}$ for $1\leq k\leq n-1$ in a fixed ascending basis.
However, by nature of an ascending basis (\ref{d omega n}), these items are equal to zero. It follows that for all $q\geq 1$.
\begin{equation}\label{dbar lambda close}
\dbar_{\Lambda}(\Gamma^{p,q}- X \wedge \alpha^{p, q-1})=0.
\end{equation}
It shows that if
\[
{\Gamma}^{p,q}=\oomega^n\wedge \alpha^{p, q-1}+\beta^{p,q},
\]
represents a class in $H^q(\lieg^{p,0})$, then
\begin{equation}\label{lifting}
\hat{\Gamma}=-X\wedge \alpha^{p, q-1}+\Gamma^{p,q}
\end{equation}
 is $\dbar_{\Lambda}$-closed. It represents a class
in
\[
\frac{F^{k-q}H^{p+q}}{F^{k-q+1}H^{p+q}}=\frac{F^{p}H^{p+q}}{F^{p+1}H^{p+q}}
\]
such that $\psi(\hat{\Gamma})=\Gamma^{p,q}$.
Therefore the map $\psi$ is an isomorphism. The map $\phi$ is its inverse. 
\eproof

\

\begin{proposition}\label{characterization}  A vector $X$ is a potential vector  of the canonical Poisson structure $V_n\wedge V_{n-1}$
associated to an ascending basis
$\{\oomega^1, \dots, \oomega^n\}$ if and only if 
\begin{equation}
\lbra{X}{\omega^k}=\iota_Xd\omega^k=0 \mbox{ for all } 1\leq k\leq n-1, \quad {\mbox { and } } \quad  \iota_Xd\omega^n=-\iota_{V_{n-1}}d\oomega^n.
\end{equation}
\end{proposition}
\bproof
Identity (\ref{dbar X}) states that
\[
\dbar X=-\sum_{k=1}^n V_k\wedge (\iota_Xd\omega^k).
\]
Since
\begin{equation}
\adL(\oomega^n)=V_n\wedge \lbra{V_{n-1}}{\oomega^{n}}=V_n\wedge \iota_{V_{n-1}}d\oomega^{n}
\end{equation}
and $\{V_1, \dots, V_n\}$ is a basis for $\lieg^{1,0}$, the conclusion of the proposition follows.
\eproof

\noindent{\textbf{Remark} }\  The above observation corrects an error in the statement of Proposition 9.2 in \cite{PS} for omitting
a requirement on the step of the nilmanifold.

\section{Identification of Gerstenhaber Algebras}

 Let
 \[
\Gamma=(\Gamma^{p+q,0}, \Gamma^{p+q-1,1}, \dots, \Gamma^{p+1, q-1}, \Gamma^{p, q}, \dots, \Gamma^{1, p+q-1},  \Gamma^{0, p+q})
 \]
 represents an element in $ \oplus_{p+q=k}H^q(\lieg^{p,0})$, the proof of Proposition \ref{surjective} shows that it is lifted to
\begin{eqnarray}
&&\phi(\Gamma^{p+q}) \nonumber \\
&=&\phi(\Gamma^{p+q,0}, \Gamma^{p+q-1,1}, \dots, \Gamma^{p+1, q-1}, \Gamma^{p, q}, \Gamma^{p-1, q+1},
\dots, \Gamma^{1, p+q-1},  \Gamma^{0, p+q})
\nonumber\\
&=&\Gamma^{p+q,0}+(\Gamma^{p+q-1, 1}-X\wedge \alpha^{p+q-1, 0})+ \dots \nonumber \\
&\  & + (\Gamma^{p, q}-X\wedge \alpha^{p, q-1})  \dots  \nonumber \\
&\  &+(\Gamma^{1, p+q-1}-X\wedge \alpha^{1, p+q-2})+ ( \Gamma^{0, p+q}-X\wedge\alpha^{0, p+q-1}).
\end{eqnarray}
Since $\Gamma^{p+q, 0}$ is in $\lieg^{{p+q,0}}$ and the complex structure is abelian, it is $\dbar_{\Lambda}$-closed so long as
it is $\dbar$-closed. With Identity (\ref{dbar lambda close}) above, we conclude that the map $\phi$ takes $\Gamma^{p+q}$ into
$H^{p+q}_\Lambda$.

If one arranges $\phi(\Gamma^{p+q})$ in terms of the bi-degree $(p, q)$ in the decomposition of
$\oplus_{p+q}B^{p,q}$, then
\begin{eqnarray}
&&\phi(\Gamma^{p+q}) \nonumber \\
&=&(\Gamma^{p+q,0}-X\wedge \alpha^{p+q-1, 0})+ (\Gamma^{p+q-1,1}-X\wedge \alpha^{p+q-2, 1}) \dots \nonumber \\
&\  & +(\Gamma^{p, q}-X\wedge \alpha^{p-1, q})+ \dots  \nonumber \\
&&
+(\Gamma^{1, p+q-1}-X\wedge \alpha^{0, p+q-1})+  \Gamma^{0, p+q}.
\end{eqnarray}

Suppose that $\adL(\oomega^n)$ is $\dbar$-exact and $X$ is a potential vector, so that
$\phi$ is an isomorphism of vector spaces from $\oplus_{p+q=k}H^q(\lieg^{p,0})$ to $H^k_\Lambda$. Recall that
$\oplus_{p+q=k}H^q(\lieg^{p,0})$ inherits the structure of a Gerstenhaber algebra. In fact, one may consider it to be the Gerstenhaber algebra
of a trivial (zero) holomorphic Poisson structure; $H^k_\Lambda$ has its own Gerstenhaber algebra structure. In the next theorem, we identify a
situation when $\phi$ is an isomorphism of Gerstenhaber algebras.

\begin{theorem}\label{algebra iso} Suppose that $\Lambda$ is the canonical holomorphic Poisson
structure associated to an ascending basis of an abelian
complex structure on a nilmanifold. If $\adL(\oomega^n)=\dbar X$ for a vector $X$,
then the map $\phi$ is an isomorphism of Gerstenhaber algebras if the action of $\ad_X=\lbra{X}{-}$ on $\liet^{*(0,1)}$ is identically zero.
\end{theorem}
\bproof
Recall  the decomposition
\begin{equation}
B^{p,q}=\{\oomega^n\}\otimes \liet^{*(p, q-1)}\oplus \liet^{p,0}\otimes\liet^{*(0,q)}.
\end{equation}
When $\Gamma$ is in $B^{p,q}$ and $\Upsilon$ is in $B^{k, \ell}$, then they have decompositions
\begin{equation}\label{decomposition}
\Gamma=\oomega^n\wedge\alpha+\beta, \quad \mbox{ and } \quad \Upsilon=\oomega^n\wedge\epsilon+\gamma,
\end{equation}
where $\alpha\in \liet^{*(p,q-1)}$, $\beta\in \liet^{*(p,q)}$, $\epsilon\in \liet^{*(k,\ell-1)}$ and 
$\gamma\in \liet^{*(k,\ell-1)}$.  
Their exterior product has a corresponding decomposition.
\[
\Gamma\wedge\Upsilon= \oomega^n\wedge(\alpha\wedge \gamma+(-1)^{p+q}\beta\wedge\epsilon)
+\beta\wedge \gamma.
\]
It follows (\ref{lifting}) that
\[
\phi(\Gamma\wedge\Upsilon)=\Gamma\wedge \Upsilon-X\wedge (\alpha\wedge \gamma+(-1)^{p+q}\beta\wedge\epsilon).
\]
On the other hand, (\ref{lifting}) and (\ref{decomposition}) together implies that 
\begin{equation}\label{lifting 2}
\phi(\Gamma)=\Gamma-X\wedge \alpha= \oomega^n\wedge \alpha+\beta-X\wedge\alpha
\end{equation}
and
\begin{equation}\label{lifting 3}
\phi(\Upsilon)=\Upsilon-X\wedge \epsilon=
\oomega^n\wedge \epsilon +\gamma-X\wedge\epsilon.
\end{equation}
It follows that 
\begin{eqnarray*}
&& \phi(\Gamma)\wedge \phi(\Upsilon)\\
&=&\Gamma\wedge\Upsilon-(\oomega^n\wedge \alpha+\beta)\wedge X\wedge\epsilon
- X\wedge\alpha\wedge (\oomega^n\wedge \epsilon +\gamma)\\
&=&\Gamma\wedge \Upsilon-X\wedge (\alpha\wedge \gamma+(-1)^{p+q}\beta\wedge\epsilon)\\
&=&\phi(\Gamma\wedge \Upsilon).
\end{eqnarray*}

Next we compare $\phi(\lbra{\Gamma}{\Upsilon})$ with $\lbra{\phi(\Gamma)}{\phi(\Upsilon)}$.
Making use of (\ref{lifting 2}) and (\ref{lifting 3}) and after  a tedious computation, we verify that
\begin{eqnarray}
&& \lbra{\phi(\Gamma)}{\phi(\Upsilon)}-\phi(\lbra{\Gamma}{\Upsilon})\nonumber \\
&=& -\oomega^n\wedge\lbra{\alpha}{X}\wedge\epsilon-\lbra{\beta}{X}\wedge\epsilon
+X\wedge\lbra{\alpha}{X}\wedge\epsilon \nonumber \\
&&-(-1)^{p+q+k+\ell}\alpha\wedge\lbra{X}{\epsilon}\wedge\oomega^n
 +(-1)^{p+q}\alpha\wedge\lbra{X}{\gamma} \nonumber\\
&&+(-1)^{p+q+k+\ell}\alpha\wedge\lbra{X}{\epsilon}\wedge X.
\end{eqnarray}

Since the complex structure is abelian, the action of $\ad_X=\lbra{X}{-}$ on $\liet^{p,0}$ is identically zero for all $p$. Therefore, if
the action of $\ad_X$ on $\liet^{*(0,1)}$ is also identically zero, then $\ad_X$ is identically zero on $\liet^{p,0}\otimes\liet^{*(0,q)}$ for all $p, q\geq 0$.
It follows that
\[
\lbra{\phi(\Gamma)}{\phi(\Upsilon)}=\phi(\lbra{\Gamma}{\Upsilon}).
\]
\eproof

\section{Applications}

\subsection{Two-steps case}\label{two-steps case}

As an application of various observations in the previous sections, we focus on two-steps nilmanifold with abelian complex structure. 
Suppose that the map 
 \[
d\omega^n: \liet^{1,0} \to \liet^{*(0,1)}
\]
 is non-degenerate, so that  $d\oomega^n$ is non-degenerate as well. It follows that $-\iota_{V_{n-1}}d\oomega^n$ is a non-zero vector. As $d\omega^n$ is non-degenerate, there 
exists $X$ in $\liet^{1,0}$ such that 
\[
\iota_Xd\omega^n=-\iota_{V_{n-1}}d\oomega^n.
\]
By nature of an ascending basis and the assumption that the algebra is two-step, $d\omega^k=0$ for $1\leq k\leq n-1$. 
It follows that all conditions in Proposition \ref{characterization} are satisfied. Therefore, $X$ is a potential vector. 

Furthermore, as $d\oomega^k=0$ for all $1\leq k\leq n-1$ as well, 
it follows that $\lbra{X}{\oomega^k}=\iota_X d\oomega^k=0$ for all 
$\oomega^k$ with $1\leq k\leq n-1$. As they span $\liet^{*(0,1)}$, then the potential vector satisfies the conditions in Theorem \ref{algebra iso}.
As a consequence, we obtain the following observation. 

\begin{theorem}\label{key} On a two-step nilmanifold with abelian complex structure, if $d\omega^n$ is non-degenerate, 
the canonical Poisson cohomology with respect to an
ascending basis is isomorphic as Gerstenhaber algebras to the one associated to the trivial (zero) Poisson structure.  i.e.
\[
\left( H^k_{\Lambda}(M), \lbra{-}{-}, \wedge \right) \cong \left(\oplus_{p+q=k}H^q(M, \Theta^p), \lbra{-}{-}, \wedge \right).
\]
\end{theorem}

\

The theorem above provides a theoretical explanation on a collection of examples in \cite{PS}. 

\subsection{Four-dimension example}
The only non-trivial real four-dimension example is the Kodaira surface, whose corresponding Lie algebra is a two-step nilpotent algebra \cite{Salamon}.
In terms of an ascending basis $\{\omega^1, \omega^2\}$, its structure equations are
\[
d\omega^1=0, \qquad d\omega^2=C\omega^1\wedge \oomega^1
\]
where $C$ is a non-zero constant. When $\{V_1, V_2\}$ is the dual basis, it is apparent that
\[
\iota_{V_1}d\oomega^2=-{\overline{C}}\oomega^1, \quad
\iota_{V_1}d\omega^2=C\oomega^1.
\]
Then $X=-\frac{\overline{C}}{C}V_2$ is a potential vector. In this case, the isomorphism is 
at the level of Gerstenhaber algebras as predicted in Theorem \ref{key} above, which provides a theoretical 
explanation for a computation in \cite{Poon}.

\subsection{A counterexample}\label{counter example}
Here we present an example to illustrate that the issue in this paper is non-trivial.

Consider a real vector space $W_6$ spanned by
$\{X_1, X_2, X_3, X_4, Z_1, Z_2\}.$ Define a Lie bracket by
\[
\lbra{X_1}{X_3}=-\frac12 Z_1, \quad \lbra{X_1}{X_4}=-\frac12Z_2, \quad \lbra{X_2}{X_3}=-\frac12 Z_2, \quad \lbra{X_2}{X_4}=\frac12 Z_1.
\]
Consider $J\circ J=-$identity and
\[
JX_1=X_2, \quad JX_4=X_3, \quad JZ_2=Z_1.
\]
It is an abelian complex structure. A basis for $\lieg^{1,0}$ is
\begin{equation}
V_1=\frac12(X_1-iX_2), \quad V_2=\frac12(X_3+iX_4), \quad V_3=\frac12(Z_1+iZ_2).
\end{equation}
The corresponding complex structure equations become
\[
\lbra{\overline{V}_1}{V_2}=-\frac12V_3, \qquad \lbra{\overline{V}_2}{V_1}=\frac12 {\overline{V}}_3.
\]
Let $\{\omega^1, \omega^2, \omega^3\}$ be the dual basis for $\lieg^{*(1,0)}$, then
\begin{equation}
d\omega^1=0, \quad d\omega^2=0, \quad \mbox{ and } \quad d\omega^3=-\frac12\omega^2\wedge\oomega^1.
\end{equation}
Since $d\oomega^3=\frac12\omega^1\wedge\oomega^2$,
\[
\iota_{V_2}d\oomega^3=\frac12\oomega^2.
\]
The equation
\[
\iota_Xd\omega^3=\iota_{V_2}d\oomega^3
\]
does not have a solution. Therefore, the canonical holomorphic Poisson structure associated to nilmanifold with the given abelian complex
structure fails to have a decomposition as noted in Theorem \ref{degeneracy}. 

\subsection{Six-dimensional, Type I Example}
Real six-dimensional nilmanifolds with abelian complex structure could be divided into two types \cite{Ugarte}.
The structure equations of an ascending basis for type I abelian complex structures are given by
\begin{equation}\label{type I}
d\omega^1=0, \quad d\omega^2=\omega^{1\overline{1}}, \quad
d\omega^3=B\omega^{1\overline{2}}+C\omega^{2\overline{1}},
\end{equation}
where $\omega^{k{\overline\ell}}=\omega^k\wedge\oomega^\ell$. It follows that
\[
d\oomega^3=-{\overline{B}}\omega^{2\overline{1}}-{\overline{C}}\omega^{1\overline{2}}.
\]
Suppose that $X=c_1V_1+c_2V_2$ is a potential vector field of the canonical Poisson structure, it satisfies the following set of constraints.
\[
\iota_Xd\omega^1=0, \quad \iota_Xd\omega^2=0, \quad \iota_Xd\omega^3=-\iota_{V_2}d\oomega^3.
\]
It terms of coordinates, these constraints become
\begin{equation}\label{type 1 in coordinates}
c_1\oomega^1=0, \quad c_1B\oomega^2+c_2C\oomega^1={\overline{B}}\oomega^1.
\end{equation}
Therefore, $c_1=0$ and $c_2={\overline{B}}/C$ if $BC\neq 0$. In other words, if the map $d\omega^3$ from $\liet^{1,0}$ to
$\liet^{*(0,1)}$ is non-degenerate, then the canonical Poisson structure has a potential vector. Moreover, the Hodge-type decomposition
is also an isomorphism of Gerstenhaber algebras.

There are two cases when $d\omega^3$ degenerates, but is non-trivial.

If $B=0$, $\ad_{V_2}\oomega^3=0$. It follows that
\[
\dbar_{\Lambda}=\dbar+\adL=\dbar
\]
when it acts on $\oplus_{p,q}B^{p,q}$. Therefore, $\oplus_{p+q=k}H^q(\lieg^{p,0})$ is naturally identical to $H^k_\Lambda$.

If $C=0$ and $B\neq 0$, the system (\ref{type 1 in coordinates}) does not have a solution.

This set of examples extends our scope beyond two-step nilmanifolds as seen in Section \ref{counter example} above. 

\subsection{Six-dimensional, Type II Example}
The structure equations of an ascending basis for Type II abelian complex structure in six-dimension are given by
\begin{equation}\label{type II}
d\omega^1=0, \quad d\omega^2=0, \quad
d\omega^3=A\omega^{1\overline{1}}+B\omega^{1\overline{2}}+C\omega^{2\overline{1}}+D\omega^{2\overline{2}}.
\end{equation}
This is a two-step nilmanifold. Note that
\[
d\oomega^3=-{\overline{A}}\omega^{1\overline{1}}-{\overline{B}}\omega^{2\overline{1}}
-{\overline C}\omega^{1\overline{2}}-{\overline D}\omega^{2\overline{2}},
\quad \iota_{V_2}d\oomega^3=-{\overline{B}}\oomega^1-{\overline D}\oomega^2.
\]
If $X=c_1V_1+c_2V_2$ is a potential vector field, then the sole constraint is $\iota_Xd\omega^3=-\iota_{V_2}d\oomega^3$,
i.e.
\[
c_1A\oomega^1+c_1B\oomega^2+c_2C\oomega^1+c_2D\oomega^2=-{\overline{B}}\oomega^1-{\overline D}\oomega^2.
\]
Equivalently,
\begin{equation}
(c_1, c_2)
\left(
\begin{array}{cc}
A  & B \\
C & D
\end{array}
\right)
=(-{\overline{B}}, -{\overline{D}}).
\end{equation}
Therefore, if $d\omega^3$ is non-degenerate,  the canonical Poisson structure has a potential vector.

If $d\omega^3$ is a rank-1 map, there are various cases. Instead of investigating every case, we consider a
particular situation.
If $C=0$ and $D=0$, the constraints on a potential vector are reduced to
\[
c_1A=-{\overline B}, \quad c_2B=0.
\]
\begin{itemize}
\item If $A\neq 0$ and $B\neq 0$, then there is a solution: $(c_1, c_2)=(-{\overline B}/A, 0)$.
\item If $A=0$ and $B\neq 0$, then there is no  solution. This case corresponds exactly to the example found in Section \ref{counter example}.
\item If $A\neq 0$ and $B=0$, then there is a solution. $(c_1, c_2)=(0, 1)$.
\end{itemize}

\

The above observation points out the complication of the existence of potential vector when $d\omega^n$ degenerates.
It will be a subject for future analysis. 

\

\noindent{\bf Acknowledgment.}
 Y.~S. Poon thanks the organizers
for the 5th Workshop
"Complex Geometry and Lie Groups" for putting together a very interesting and rich meeting at University of Florence in spring 2018.
He thanks Anna Fino for helpful suggestions during this conference. 
He also thanks the Institute of Mathematical Sciences of the Chinese University of Hong Kong for
providing an excellent working environment in summer 2018 for him to complete this manuscript.

\noindent{Yat Sun  Poon:
    Department of Mathematics, University of California,
    Riverside, CA 92521, U.S.A.. E-mail: ypoon@ucr.edu.}

\

\noindent{John Simanyi: Department of Mathematics, University of California,
    Riverside, CA 92521, U.S.A. }

\end{document}